\documentclass[12pt,a4paper]{article}
\usepackage{theorem}
\newtheorem{lemma}{Lemma}

\newtheorem{theorem}[lemma]{Theorem}

{\theorembodyfont{\upshape}}
{\theorembodyfont{\upshape}}
{\theorembodyfont{\upshape}\newtheorem{example}[lemma]{Example}}
{\theorembodyfont{\upshape}\newtheorem{stage}{Stage}}

\newcommand{\Z}{{\bf Z}}
\newcommand{\R}{{\bf R}}
\newcommand{\C}{{\bf C}}

\newcommand{\cC}{{\cal C}}
\newcommand{\cD}{{\cal D}}
\newcommand{\cE}{{\cal E}}
\newcommand{\cF}{{\cal F}}
\newcommand{\cG}{{\cal G}}
\newcommand{\cH}{{\cal H}}

\newcommand{\cQ}{{\cal Q}}

\newcommand{\alp}{\alpha}
\newcommand{\bet}{\beta}
\newcommand{\gam}{\gamma}
\newcommand{\lam}{\lambda}
\newcommand{\del}{\delta}
\newcommand{\Gam}{\Gamma}
\newcommand{\Ome}{\Omega}
\newcommand{\lap}{{\Delta}}
\newcommand{\degree }{{\rm deg}}

\newcommand{\norm}{\Vert}

\newcommand{\Proof}{\underbar{Proof}{\hskip 0.1in}}
\newcommand{\hash}{\#}
\newcommand{\la}{{\langle}}
\newcommand{\ra}{{\rangle}}
\newcommand{\rmd}{{\rm d}}
\title{A HIERARCHICAL METHOD\break FOR OBTAINING\break  EIGENVALUE ENCLOSURES}
\author{E.B. Davies}
\date{ January 1998}
\begin{document}
\maketitle
\begin{abstract}
We introduce a new method of obtaining guaranteed enclosures 
of the eigenvalues of a variety of self-adjoint 
differential and difference operators 
with discrete spectrum. The method is based upon subdividing 
the region into a number of simpler regions for which 
eigenvalue enclosures are already available.
\end{abstract}
\section{Introduction}
\par
A rigorous method of obtaining enclosures of the 
eigenvalues of self-adjoint operators 
has recently been described by Goerisch and Plum 
\cite{G, P1, P2, P3}. It depends upon having a 
soluble comparison operator, from which a controlled 
homotopy is carried out. In this paper we introduce a new 
method which has the advantage of not requiring such a 
comparison operator, and apply it to a variety of examples. 
In Sections 2 to 5 we consider 
Sturm-Liouville operators in some detail. Sections 6 and 7 describe how to 
adapt the method to higher order operators and systems, 
still in one dimension. In Sections 8 and 9 we treat 
discrete Laplacians on graphs, while in Section 10 we 
consider the Laplacian acting in a bounded region in 
Euclidean space. The method can be applied to second and 
higher order elliptic differential 
operators with variable coefficients, but we do not present 
the details here. 
\par
We distinguish between computing an eigenvalue in floating 
point arithmetic, and obtaining guaranteed enclosures. When 
using the word `enclosure' we shall always understand 
that the calculation is mathematically 
rigorous, and that the computations are done in 
interval arithmetic. Most numerical computations do not give 
proofs that the values obtained are correct, but depend upon 
the experience of the person who writes or uses the program 
concerning its range of reliability. With guaranteed 
enclosures on the other hand, the value obtained is known 
to be correct within the stated error bounds, unless there 
is an actual error at some stage of the computation.
\par 
There are already several methods of computing the 
eigenvalues of a Sturm-Liouville operator $H$ acting in 
$L^{2}(\alp,\bet)$, and higher order analogues. The most obvious 
one, called shooting, solves the initial value problem for 
the eigenvalue equation $Hf=\lam f$ subject to the given 
boundary conditions at $\alp$ and then varies $\lam$ until 
the boundary condition at $\bet$ is also valid. Most 
shooting programs do not try to give guaranteed error 
bounds. Although this is entirely possible \cite{Lo}, the 
method is difficult to implement numerically if 
the potential is singular at both ends of the interval. 
\par
A second method, introduced by Goerisch \cite{G} and 
Plum \cite{P1,P2,P3}, obtains guaranteed 
enclosures on the eigenvalues of a self-adjoint operator 
$H$ by a continuous homotopy method, starting from a 
simpler operator. This is often exactly soluble, but a 
minimum requirement is that one can obtain sufficiently good 
rigorous lower bounds on its eigenvalues. Our method is 
similar to theirs in that it 
also uses a homotopy from a simpler operator. However they 
consider a continuous homotopy in some parameter, 
which often changes the coefficients smoothly to those of an operator 
with constant coefficients, 
while we consider a discrete homotopy in certain internal boundary 
conditions which we choose to insert. 
We have compared our variation of the homotopy method with 
theirs for some of the examples Plum solves, and it appears 
to be substantially more efficient. In higher dimensions we are 
able to treat examples which are beyond the earlier 
method, because of the non-existence of an exactly 
soluble operator possessing a continuous homotopy to the given operator.
\par
One may obtain rigorous upper bounds on any specified 
number of eigenvalues by means of the 
Rayleigh-Ritz (RR) or variational method \cite{D1}. The starting point is the 
determination of accurate approximations to the eigenfunctions 
by a non-rigorous auxiliary calculation, possibly an inverse iteration 
method. Once these have been found one starts again using RR 
to obtain rigorous upper bounds on the 
eigenvalues of the operator $H$ in interval arithmetic. 
\par
The lower bound is obtained by the method of Temple-Lehmann 
(TL) which also depends upon the choice of suitable test 
functions \cite{D1,D2,P3,ZM}. However, in this case one also needs to have crude 
lower bounds on the eigenvalues, and these are precisely 
what is missing at the rigorous level. More precisely if 
the eigenvalues of $H$, written in increasing order and 
repeated according to multiplicity, are 
$\{ \lam_{n} \}_{n=0}^{\infty}$, then in order to obtain an 
accurate lower bound on $\lam_{n}$ for some $n$ using TL one needs already to be 
in possession of a number $\rho$ such that 
\[
\lam_{n} <\rho <\lam_{n+1}
\]
where $\rho$ is not too close to $\lam_{n}$. There are 
three possible methods of obtaining such lower bounds.
\par
(i) One might hope 
that the upper bound on $\lam_{n+1}$ is fairly accurate and take 
$\rho$ to be a slightly smaller number. This idea cannot be 
turned into a rigorous procedure and will not be discussed 
further.
\par
(ii) One can use the Goerisch-Plum coefficient homotopy 
method of obtaining enclosures for many operators.
\par
(iii) One can use a boundary condition homotopy method. The 
description of this new method is the main contribution of this paper.
\par
In the above we have not mentioned the extra complications 
which arise if $\lam_{n}$ is degenerate or nearly so. There 
are well-known modifications of TL which deal with this 
problem \cite{D1,P1,ZM}, but we did not want to over-complicate the 
discussion at this stage.
\par
\section{Neumann decoupling}
\par
Let $H$ be a Sturm-Liouville operator acting 
in $L^{2}(\alp,\bet)$. We 
assume that $H$ is of the form
\[
Hf(x):=-{\rmd\over\rmd x}\left\{ a(x){\rmd f\over \rmd 
x}\right\} +V(x)f(x)
\]
where $a$ is a positive function in $C^{1}[\alp,\bet]$ 
and $V\in L^{\infty}[\alp,\bet]$. We 
assume Neumann boundary conditions (NBC) in order to 
emphasise that the method does not depend upon the very 
strong monotonicity properties which hold for Dirichlet 
boundary conditions \cite{D1}. 
\par
Our method is based upon decoupling the interval 
$(\alp,\bet)$ into $2^{N}$ subintervals; in most examples 
considered by the author one can take $N=3$ or $N=4$. The 
subintervals do not need to be of equal length, but this is 
the easiest choice to make. We put 
\[
\alp=\alp_{0}<\alp_{1}<\ldots<\alp_{2^{N}}=\bet
\]
and let $H_{i}$ denote the restriction of $H$ to 
$L^{2}(\alp_{i-1},\alp_{i})$ subject to NBC. We then define 
$A_{N}$ to be the sum of the $H_{i}$, so that $A_{N}$ once 
again acts in $L^{2}(\alp,\bet)$. The operators $H$ and 
$A_{N}$ have the same quadratic form 
\[
Q(f):=\int_{\alp}^{\bet}\{ 
a(x)|f^{\prime}(x)|^{2}+V(x)|f(x)|^{2}\} \rmd x
\]
but with different quadratic form domains $\cQ(H)$ and 
$\cQ(A_{N})$. The space of 
all $C^{1}$ functions on $[\alp,\bet]$ is a quadratic form 
core for $H$, but to obtain a quadratic form core for $A_{N}$ 
one must allow the functions to have arbitrary jump 
discontinuities at each $\alp_{i}$. Since  
$\cQ(H)\subset\cQ(A_{N})$, the RR method \cite{D1} shows that the 
eigenvalues of $A_{N}$ are less than or equal to those of $H$.
\par
We define intermediate operators $A_{n}$ acting in 
$L^{2}(\alp,\bet)$ for $0\leq n\leq 
N$, by a similar method. The component operators 
of $A_{n}$ are similar to those of $A_{N}$ but 
only using the points $\alp_{i}$ where $i=j.2^{N-r}$ for $0\leq j\leq 
2^{r}$. At each stage the quadratic form domain decreases, 
leading to the operator inequalities
\[
A_{N}\leq A_{N-1}\leq \ldots\leq A_{0}=H.
\] 
The following theorem enables the eigenvalues of $A_{n}$ 
to be computed rigorously using TL once one knows those of 
$A_{n+1}$. Since the passage from $A_{n+1}$ to $A_{n}$ 
consists of putting together intervals in pairs, and the set of  
eigenvalues of $A_{n}$ is simply the collection of all 
eigenvalues of its component operators $H_{i}$, it is 
sufficient to deal with the following special case, which 
also describes the passage from $A_{1}$ to $A_{0}=H$.
\par
Let $\alp<\gam<\bet$ and let $\{ \lam_{i}\}$, $\{ \mu_{i}\}$, 
$\{ \nu_{i}\}$, denote the eigenvalues of the operators 
$H_{1}$, $H_{2}$ and $H$ 
associated with the intervals $(\alp,\gam)$, $(\gam,\bet)$, 
$(\alp,\bet)$ respectively, all subject to NBC. Finally 
let 
\[
\{ \sigma_{i}\}:=\{ \lam_{i}\}\cup\{ \mu_{i}\}
\] 
subject to re-ordering in increasing order and repeating 
according to multiplicities. Then $\{ \sigma_{i}\}$ are the 
eigenvalues of the operator $A:=H_{1}+H_{2}$.
\par
\begin{theorem}
The eigenvalues $\{ \sigma_{i}\}$ and $\{ \nu_{i}\}$ 
interlace in the sense that 
\[
\sigma_{i}\leq \nu_{i}\leq \sigma_{i+1}
\]
for all $i$. Moreover these are strict 
inequalities unless the derivative of the relevant 
eigenfunction of $H$ vanishes at the point $\gam$.
\end{theorem}
\Proof  The idea is that $H$ differs from $A$ by a rank 
one perturbation in a certain singular sense. More precisely let 
$s<\sigma_{0}$ and compare $(A+s)^{-1}$ with $(H+s)^{-1}$. 
Both have Green functions which can be computed from the 
two fundamental solutions of the differential equation
\[
-{\rmd\over\rmd x}\left\{ a(x){\rmd f\over \rmd x}\right\} 
+V(x)f(x)+sf(x)=0.
\]
If one computes the difference of the two kernels one finds 
that it is a rank one operator. The inequality which we 
want is equivalent to 
\[
(\sigma_{i}+s)^{-1}\geq (\nu_{i}+s)^{-1}\geq 
(\sigma_{i+1}+s)^{-1}
\]
and this holds whenever one has a positive rank one perturbation, by 
an application of the min-max principle. For an alternative proof 
see Theorem 5.
\par
The eigenvalues of $H_{1}$ are all distinct, as are the 
eigenvalues of $H_{2}$, because $H_{i}$ are Sturm-Liouville 
operators. However there may be coincidences between the two 
sets of eigenvalues, which imply that 
$\sigma_{i}=\sigma_{i+1}$. This is not a serious problem, 
but it can usually be avoided if desired by moving 
the point $\gam$ slightly. Assuming that this has been done 
it is not possible that $\sigma_{i}= \nu_{i}$: this 
equality would imply that 
the corresponding eigenfunction of $H$ happens to have zero 
derivative at $\gam$, in which case it is also an 
eigenfunction for both $H_{1}$ and $H_{2}$, so $\sigma_{i}=\sigma_{i+1}$. 
\par 
\section{The Enclosure Algorithm}
\par
We start with a subdivision of $(\alp,\bet)$ into $2^{N}$ 
parts which is fine enough for us to be able to obtain 
disjoint enclosures on 
the eigenvalues of each component $H_{i}$ by comparison with constant 
coefficient operators. From here onwards we suppose that we 
are only interested in obtaining enclosures on those 
eigenvalues of $H$ which are less than some pre-assigned number 
$E$. The larger the value of $E$, the larger one must take 
$N$ in order to be able to start the procedure. The following 
lemma shows that it is 
sufficient for the coefficients to be close to constant in 
each interval. We only consider the case of $H$ itself for 
notational simplicity, but the lemma should actually be 
applied to each component $H_{i}$ of $A_{N}$.
\par
\begin{lemma}
Suppose that $a_{0}\leq a(x)\leq a_{1}$ and 
$v_{0}\leq V(x)\leq v_{1}$ for all $x\in 
(\alp,\bet)$. Then the eigenvalues $\{ \lam_{i}\}$ of $H$ 
satisfy
\[
a_{0}\pi^{2}i^{2}/(\bet-\alp)^{2}+v_{0 }\leq \lam_{i}
\leq a_{1}\pi^{2}i^{2}/(\bet-\alp)^{2}+v_{1 }. 
\]
These enclosure intervals are disjoint for all eigenvalues 
less than a given number $E^{\prime}$ if 
\[
0\leq v_{1}-v_{0}\leq a_{0}\pi^{2}/(\bet-\alp)^{2}
\]
and  
\[
0\leq v_{1}-v_{0}\leq M^{2}a_{0}\pi^{2}/(\bet-\alp)^{2}
-(M-1)^{2}a_{1}\pi^{2}/(\bet-\alp)^{2}
\]
where $M$ is the smallest integer such that
\[
E^{\prime}\leq v_{0}+M^{2}a_{0}\pi^{2}/(\bet-\alp)^{2}.
\]
\end{lemma}
\par
\Proof  The first inequality follows by comparing $H$ with the 
obvious constant coefficient operators, whose eigenvalues 
are exactly computable. The proof of the second uses the 
observation that
\[
v_{0}+m^{2}a_{0}\pi^{2}/(\bet-\alp)^{2}
-v_{1}-(m-1)^{2}a_{1}\pi^{2}/(\bet-\alp)^{2}
\]
is a concave function of $m$ which must therefore take its 
minimum value on the interval $[1,M]$ at one of its ends.
\par
\underbar{Note}{\hskip 0.07in} Since the initial intervals $(\alp_{i-1},\alp_{i})$ are 
quite short, the integer $M$ in the above lemma may be quite 
small and the above lemma may not impose strong 
conditions on the constants $v_{0},v_{1},a_{0},a_{1}$.
\par
The algorithm for obtaining enclosures of the eigenvalues 
of $H$ has several stages:
\par
\begin{stage} We have to choose an initial subdivision of 
the interval $(\alp,\bet)$ such that each of the 
subintervals $(\alp_{i-1},\alp_{i})$ satisfies the conditions 
of Lemma 2. This can be done in several ways and is discussed 
further in Section 4.
\end{stage}
\begin{stage} We choose a number $E^{\prime}>E$, for example 
$E^{\prime}:=9E/8$, and put $E_{n}:=E+n(E^{\prime}-E)/N$ 
for all $0\leq n\leq N+1$.
\end{stage}
\begin{stage} We subdivide $(\alp,\bet)$ as described above for a 
value of $N$ which is large enough for us to 
obtain disjoint intervals which enclose each of the 
eigenvalues of each component $H_{i}$ of $A_{N}$ up the number 
$E_{N+1}$.
\end{stage}
\begin{stage} We use RRTL to obtain accurate enclosures of each of the 
eigenvalues of each component $H_{i}$ up to the number 
$E_{N}$. Putting these 
together in pairs we obtain rough enclosures of the 
eigenvalues of each component $H_{j}$ of $A_{N-1}$ by virtue 
of Theorem 1. These 
enclosure intervals overlap very slightly because the 
previous accurate enclosures were not perfect. 
\end{stage}
\begin{stage} We apply RR to each component $H_{j}$ of 
$A_{N-1}$ to obtain smaller upper 
bounds on each of the eigenvalues of each $H_{j}$ and so to 
convert the above into 
rough but nevertheless disjoint enclosures of the eigenvalues of each 
$H_{j}$ up to $E_{N}$. 
\end{stage}
\begin{stage} We apply TL to obtain accurate enclosures of 
each of the eigenvalues of each component $H_{j}$ of 
$A_{N-1}$ up to $E_{N-1}$.
\end{stage}
\begin{stage} We repeat the process inductively until we reach 
accurate enclosures of each of the eigenvalues of $H$ up to 
$E$.
\end{stage} 
Some comments are in order.
\par
The introduction of the 
sequence $E_{n}$ at Stage 1 is needed because TL requires a 
significant gap above any eigenvalue to be estimated. If 
there is an eigenvalue very close to the upper limit 
$E_{n}$ at any stage then that eigenvalue cannot be estimated accurately.
\par
When the eigenvalues af two adjacent operators are combined 
in Stage 3 it may happen that two eigenvalues of the new 
list created coincide to a high degree of accuracy. This is one 
possible cause of the problem mention in the next paragraph.
\par
The procedure in Stage 5 may occasionally fail because RR 
may not decrease the upper bound on an eigenvalue enough to 
make the intervals disjoint. This is handled by using a 
higher order version of TL whenever this occurs. In 
principle this could occur for all eigenvalues, in which 
case the algorithm might halt, but this is extremely 
unlikely unless there is a symmetry of the underlying 
problem, which should have been taken into account before 
starting the computation. 
\par
Ultimately we do not guarantee 
either that the algorithm finishes or that the results 
which it yields are of the desired accuracy, but only that 
if the algorithm does finish then the enclosures obtained are correct. 
If the enclosures are not sufficiently accurate, 
then one must start again with a larger test function space.
\par
Although we specified that the second order coefficients 
$a(x)$ of the 
differential operator should be $C^{1}$, there is no 
difficulty in accommodating simple jump discontinuities. 
Once one has determined the location of these points, they 
should be included in the partition $\{ 
\alp_{i}\}_{i=0}^{2^{N}}$ of the interval $(\alp,\bet)$. 
The discontinuity of $a(x)$ at a point $\gam$ imposes
an effective internal boundary 
condition on $H$ at $\gam$, which must be taken 
into account when specifying its operator domain, but has 
no effect on its quadratic form domain.
\par
One way of estimating the total computational effort is to 
count the number of distinct operators for which we have to 
compute some of the eigenvalues accurately. At the level $n$ 
this is $2^{n}$, so the total number is $2^{N+1}-1$. 
\par
Since parallel machines will become more important, it 
should be noted that the computations of the eigenvalues of 
the different operators $H_{j}$ at any particular level are 
entirely independent, and may be carried out 
simultaneously. Thus on a parallel machine the total 
computational effort is proportional to $N+1$. In all the 
examples which we have considered this means the algorithm 
has only three or four steps!
\par
Both of the above estimates of computational effort are
too pessimistic. At the 
higher levels it may be seen in the examples we analyse 
below that the number of eigenvalues 
of each operator to be computed is very small, because 
the eigenvalues are far apart. So the 
computation is much faster at the higher levels than 
indicated above, whether or not one has a parallel machine.
\par
It is clear that the same procedure may be used irrespective 
of the boundary conditions at $\alp$ and $\bet$. It may also 
be applied to potentials which are singular at the end points 
provided one has crude bounds on the eigenvalues to replace 
those of Lemma 2. Its extension to systems and to higher order differential 
operators is described in Sections 6 and 7.
\par
\section{The Subdivision of $(\alp,\bet)$}
\par
We have suggested above that the subdivision of  
$(\alp,\bet)$ should be defined by 
\[
\alp_{i}:=\alp+i(\bet-\alp)/2^{N}
\]
for $0\leq i\leq 2^{N}$ ,where the size of $N$ is 
determined by Lemma 2 as indicated 
in the algorithm. However it is possible that when 
combining two eigenvalue lists one finds that two 
eigenvalues coincide or are undesirably close. This is not 
an insuperable problem since one can use a higher order 
version of TL to obtain the required lower bounds on the 
eigenvalues. However, there is a systematic way of 
avoiding it, unless one of the eigenfunctions has 
an interval of constancy.  
\par
We first emphasise that there is no hope of obtaining 
accurate enclosures of the eigenvalues of $H$ unless there 
is some other non-rigorous method of computing the 
eigenvalues, such as unsupplemented RR or shooting, which 
in fact give good approximations to the eigenvalues and 
eigenfunctions. We use these computed eigenfunctions to 
choose the bisection point $\gam$ of $(\alp,\bet)$ as 
described below. We then 
do the same for both of the subintervals $(\alp,\gam)$ and 
$(\gam,\bet)$ and so on until we have produced a fine 
enough subdivision of $(\alp,\bet)$ according to the 
criterion of Lemma 2. If we have misled ourselves about the 
best choice of the points $\alp_{i}$ then nothing is lost, 
because we can still use the above algorithm. If however, 
the approximations to the eigenfunctions are accurate enough 
then the method we now describe will have prevented the 
problem mentioned above.
\par
Let us suppose that there are $k+1$ eigenvalues 
of $H$ less than $E$, and that the corresponding 
eigenfunctions $f_{r}$ have zero derivatives at $p(r)$ points 
for each $0\leq r\leq k$. Then we choose $\gam$ 
somewhere near the centre of $(\alp,\bet)$ but not at or 
near to any 
of the above points. Since there are $P:=p(0)+\ldots+p(k)$ 
such points altogether there exists $\gam\in (3\alp/4 
+\bet/4,\alp/4+3\bet/4)$ which is at a distance at least 
$(\bet-\alp)/4P$ from each of the points.
\par
Whether or not it is worth using this iterative procedure for 
selecting the subdivision of $(\alp,\bet)$ remains to be 
seen. In the two cases solved below, it appears that using 
a uniform subdivision is perfectly satisfactory.
\par
There is an entirely different reason for choosing a 
non-uniform subdivision of $(\alp,\bet)$, if the 
coefficients of $H$ vary substantial from one part of the interval 
to another. If the potential is bigger than the 
number $E$ in some interval, then one should make that 
entire interval one of the $(\alp_{i-1},\alp_{i})$, however 
big it is, because there will be no relevant eigenvalues 
associated with it. More generally the size of each 
interval $(\alp_{i-1},\alp_{i})$ should 
be as big as possible subject to being able to obtain 
disjoint enclosures of all of the eigenvalues of $H_{i}$ less than $E$. This 
procedure reduces the number of operators for which one has 
to compute some of the eigenvalues. A more thorough 
investigation might involve the uncertainty 
principle, but Lemma 2 suffices for most purposes.    
\section{ Examples}
\par
We illustrate our general theory with two numerical 
examples, which are solved using shooting and floating point 
arithmetic, not using interval arithmetic as is actually 
required. 
There are two reasons for this, the first being that our 
goal here is only to examine the feasibility of the method, not to 
create a new software package. The second is that 
one should not use a high-powered technique for obtaining 
eigenvalue enclosures until one has a good idea of the approximate 
location of the eigenvalues. This information cannot be 
used in the final computation because it is not rigorous, but 
it may indicate problems which need special attention in 
the rigorous computation.
\par
The two examples were studied in detail by Plum \cite{P1} 
using a continuous homotopy in the coefficients.
\par
\begin{example} Let $H$ be the operator defined by
\[
Hf(x):=-{\rmd^{2} f\over \rmd x^{2}}+8\cos(x)^{2}f(x)
\] 
acting in $L^{2}(0,\pi )$ subject to NBC. Plum obtained 
enclosures on the eigenvalues ranging from
\[
\mu_{0}=2.48604311^{50}_{47}
\]
to
\[
\mu_{8}=68.03175^{8}_{6}
\]
using his homotopy method, RRTL and interval 
arithmetic. Putting $N=2$ the conditions of Lemma 2 are 
satisfied for any choice of $E$ for each of the components 
$H_{i}$, $1\leq i\leq 4$, with $v_{1}-v_{0}= 4$, 
$a_{0}=a_{1}=1$ and $\alp_{i}-\alp_{i-1}=\pi/4$.
Now let $K_{1},K_{2}$ be the two operators at level one, 
acting in the intervals $(0,\pi/2)$ and $(\pi/2,\pi)$. We 
have computed the eigenvalues of all of the operators above 
up to the limit $E=70$.
\par
For $i=1,4$ Lemma 2 yields the crude enclosures
\[
4<\mu_{0}<8,\, 20<\mu_{1}<24,\, 68<\mu_{2}<72
\]
while more accurate, but non-rigorous, calculations provide
\[
\mu_{0}\simeq 6.454,\, \mu_{1}\simeq 22.450,\,  \mu_{2}\simeq 70.515
\]
For $i=2,3$ Lemma 2 yields the crude enclosures
\[
0<\mu_{0}<4,\,  16<\mu_{1}<20,\, 64<\mu_{2}<68,\, 144<\mu_{3}
\]
while more accurate calculations yield
\[
\mu_{0}\simeq 1.364,\, \mu_{1}\simeq 17.693, \, \mu_{2}\simeq 65.503
\]
We now join together the eigenvalue lists of $H_{1}$ and 
$H_{2}$ to obtain the list
\[
1.364,\,6.454,\,17.693,\,22.450,\,65.503,\,70.515
\] 
If these values are indeed accurate then according 
to Theorem 1 they interlace the eigenvalues of 
$K_{1}$ and provide the basis for the use of RRTL for 
obtaining accurate enclosures of the eigenvalues of 
$K_{1}$. The eigenvalues of $K_{1}$ are (again 
non-rigorously)
\[
2.486,\, 9.173,\, 20.141,\, 40.057,\, 68.032
\]
These coincide with the eigenvalues of $K_{2}$ since we 
have not made use of the symmetry of the operator about $x=\pi/2$. 
When we combine this list with a second copy of itself the 
resulting list interlaces the eigenvalues of $H$, 
namely
\[
2.486,\, 6.397,\, 9.173,\, 13.370,\, 20.141,\, 29.084,\, 
40.057,\, 53.042,\, 68.032
\]
It would have been possible to avoid the coincidence of the 
eigenvalues of $K_{1}$ and $K_{2}$ by starting with the 
partition $0,\, 0.6,\, 1.2,\, 2.1,\, \pi$ instead of the 
partition into equal subintervals.
\par
The above computation involves determining the eigenvalues 
of $7$ operators at $3$ different levels. At the top level 
we only had to compute the first $3$ eigenvalues of each 
operator $H_{i}$, at the 
middle level we had to compute $5$ and at the bottom level we 
had to compute $9$.
\par
By comparison Plum computed the eigenvalues of $10$ 
intermediate operators, with less scope for parallelization 
since each computation depended on the previous one. We 
computed a total of $31$ eigenvalues, while Plum computed at 
least $90$. Plum was not, however, particularly concerned 
with minimising numerical effort in his paper.
\end{example}
\par
\begin{example}
We consider the operator 
\[
Hf(x):=-{\rmd^{2} f\over \rmd x^{2}}+1000xf(x)
\]
acting on $L^{2}(0,1)$ subject to DBC. This is essentially 
the same as Example 2 of Plum \cite{P1}, who obtained enclosures 
on the eigenvalues ranging from
\[
\mu_{0}=233.8107^{42}_{35}
\]
to
\[
\mu_{9}=1508.10^{83}_{78}
\]
The unpublished enclosures of Lohner \cite{Lo}, obtained by shooting, 
are considerably more accurate.
We put $N:=3$, $\alp_{i}:=i/8$ for $0\leq i\leq 8$ and 
$E:=1000$. A 
modification of Lemma 2 to cope with the DBC at $0,1$ 
yields the following crude eigenvalue enclosures. 
\par
In $(0,1/8)$ we have initially
\[
157<\mu_{0}<283,\, 1421<\mu_{1}<1547
\]
and then more accurately
\[
\mu_{0}\simeq 245.225,\, \mu_{1}\simeq 1486.798
\]
In $(1/8,1/4)$ we have initially
\[
125<\mu_{0}<250,\, 756<\mu_{1}<882
\]
and then more accurately
\[
\mu_{0}\simeq 185.471,\, \mu_{1}\simeq 820.761
\]
We omit the results for the other intervals, each of which 
involves the computation of only two eigenvalues,  all 
higher eigenvalues being bigger than 1500.
We now consider level 2. Putting the previous lists together 
in pairs and considering the interval $(0,1/4)$ 
we obtain the crude bounds
\[
185.471\leq\mu_{0}\leq 245.225\leq \mu_{1}\leq 820.761 \leq 
\mu_{2 }\leq 1486.798 
\leq \mu_{3}
\]
In fact the upper bound on each $\mu_{i}$ is slightly 
bigger than the lower bound on $\mu_{i+1}$, because the 
number separating them is not exact, but the use of RR reduces the upper 
bound on each eigenvalue substantially and so yields disjoint enclosures of 
the eigenvalues (or actually would do so if the calculations 
were rigorous). All higher eigenvalues are greater than 
$1500$. The accurate eigenvalues for the interval $(0,1/4)$ 
are
\[
\mu_{0}=205.942,\, \mu_{1}=490.938,\, \mu_{2}=1115.419 
\]
We omit the further computations at levels 2 and 1. The full list 
of eigenvalues of $A_{1}$ is 
\[
233.705,\, 400.348,\,  532.152,\,  601.881,\,  748.110,\,  825.999,\, 1007.897
\]
and interlaces the list of eigenvalues of $H$, namely:
\[
233.811,\,  408.795,\,  552.056,\,  678.679,\,  794.738,\,  906.461
\]
which agree with the enclosures of Plum \cite{P1}. We next comment on 
the amount of computation needed by our method.
\par
The total number of operators considered by our method, 
is $15$, compared with $50$ in Plum's method, 
since he puts $\delta=0.02$. If one has a parallel machine 
then the relevant quantity is the number of levels, namely 
4. For each operator at level 3 we needed to compute $1$ or $2$ eigenvalues. 
For each operator at level 2 we computed $3$ eigenvalues.
For the two operators at level 1 we computed $5$ 
and $3$ eigenvalues. 
Finally we computed all $6$ eigenvalues of $H$ in the 
interval $[0,1000]$, 
making a total of at most $42$ eigenvalues computed. Plum's method 
involves the computation of at least $300$ eigenvalues, but 
he did not attempt to minimise this number.   
\end{example}
\par
\section{Higher Order Operators}
\par
The procedure which we described above can be modified to 
treat higher order differential operators in one dimension. 
The difference in the higher order case is that decoupling 
an interval into two parts by introducing a Neumann 
boundary condition is not equivalent to a rank one 
perturbation. However it is still of finite rank, as we will 
now explain. 
\par
Let $H$ be defined formally on $L^{2}(\alp,\bet)$ by
\[
Hf(x):=(-1)^{m}{\rmd^{m} \over \rmd x^{m}} \left\{ a(x) {\rmd^{m} f\over 
\rmd x^{m}}\right\}
\]
where $a\in C^{m}[\alp,\bet ]$ is positive. Our method can also 
deal with more complicated operators involving lower order 
terms. We assume 
Neumann boundary conditions, in the sense that we take the 
quadratic form of the operator to be 
\[
Q(f):=\int_{\alp}^{\bet }a(x)\left|{\rmd^{m} f\over 
\rmd x^{m}}\right|^{2}\rmd x
\]
with domain the Sobolev space $W^{m,2}(\alp,\bet)$. It is 
known that $Q$ is closed on this domain, and we define $H$ 
to be the non-negative self-adjoint operator associated with 
the form in the standard manner.
\par
Functions in  $W^{m,2}(\alp,\bet)$ are continuous on 
$[\alp,\bet]$ along with all derivatives of order less than 
$m$. Given $\alp<\gam<\bet$ we introduce a Neumann boundary 
condition at $\gam$ by replacing $W^{m,2}(\alp,\bet)$ by 
the space $\cQ_{m}$ in which we allow the functions and their first 
$m-1$ derivatives to have simple jump discontinuities at 
$\gam$. Let $H_{m}$ be the corresponding operator on  
$L^{2}(\alp,\bet)$. We now define a chain of operators 
$H_{r}$ for $0\leq r\leq m$ with $H_{0}:=H$. Each of them 
is associated with the same form $Q$ but on different 
domains $\cQ_{r}$. We define $\cQ_{r}$ to be the space of 
functions in $\cQ_{m}$ such that all derivatives of $f$ from 
the order $r$ to $m-1$ inclusive are continuous at 
$\gam$. Thus $\cQ_{r}\subset \cQ_{r+1}$ for all $r$, each 
being of co-dimension one in the next.
\par
\begin{theorem}
Let $H$ and $K$ be two non-negative self-adjoint operators  
on a Hilbert space $\cH$ such that their quadratic forms 
coincide on their common domain. Suppose also that the form 
domain $\cQ(K)$ of $K$ is a subspace of co-dimension $1$ in 
$\cQ(H)$. Finally suppose that $H$ and $K$ both have purely 
discrete spectrum  and that their eigenvalues
written in increasing order and repeated according 
to multiplicity are respectively $\{ \lam_{n}\}_{n=0}^{\infty}$ 
and  $\{ \mu_{n}\}_{n=0}^{\infty}$. 
Then the two sets of eigenvalues interlace in the sense that
\[
\lam_{n}\leq \mu_{n}\leq \lam_{n+1}
\]
for all $n$.  
\end{theorem}
\Proof
This is an immediate consequence of the min-max principle, since 
every subspace of dimension $n$ of $\cQ(H)$ is either 
already contained in $\cQ(K)$ or intersects $\cQ(K)$ in a 
subspace of dimension $n-1$. 
\par
The application of this theorem to higher order operators is 
immediate. In order to remove a Neumann boundary condition 
at the point $\gam$ we have to pass through a chain of 
operators $H_{r}$ with $r$ decreasing from $m$ to $0$. At 
each stage the eigenvalues interlace, and this is the 
condition needed to apply the TL technique as described in 
Section 3.
\par
We have described the operator $H_{r}$ in terms of its 
quadratic form domain. This is sufficient for the 
application of the RR technique. However, the TL method 
depends upon the selection of test functions from the 
operator domain, so we need to describe this. We first 
comment that functions in any of the operator domains lie in 
$C^{2m-1}[\alp,\gam]+C^{2m-1}[\gam,\bet]$ because of our smoothness assumption 
on the coefficient $a(x)$. Also the weak derivative 
$f^{(2m)}(x)$ in each subinterval 
must lie in $L^{2}$. Of course eigenfunctions are 
more regular and must lie in 
$C^{2m}[\alp,\gam]+C^{2m}[\gam,\bet]$.
\par
We now specify the boundary conditions. The choice of  
quadratic form domain implies that if $f$ lies in the operator domain then 
$f^{(r)}(x)=0$ for $x=\alp,\bet$ and for all 
$m\leq r\leq 2m-1$, i.e. Neumann boundary conditions.  We need to 
impose $2m$ boundary conditions at $\gam\pm$ to obtain a self-adjoint 
operator, and these are 
different for each operator $H_{r}$. Our quadratic form  
assumption is that $f^{(s)}(\gam-)=f^{(s)}(\gam+)$ for 
all $s$ such that $r\leq 
s \leq m-1$. This corresponds to the assumption that
\[
\begin{array}{rcll}
f^{(s)}(\gam+)&=&f^{(s)}(\gam-)&\mbox{for all $r\leq s\leq 
m-1$}\\
(af^{(m)})^{(s)}(\gam+)&=& (af^{(m)})^{(s)}(\gam-)& \mbox{ for all 
$0\leq s\leq m-r-1$}\\
(af^{(m)})^{(s)}(\gam\pm)&=& 0&\mbox{for all $m-r\leq s\leq m-1$}\\
\end{array}
\]
for all $f$ in the operator domain of $H_{r}$, as one may see 
by carrying out some integrations by parts and requiring the 
boundary terms to vanish.
\par
The test functions chosen for the TL procedure must satisfy all of the 
above boundary conditions. One could use a space 
consisting of different polynomials in each subinterval, 
with the coefficents restricted to satisfy the boundary 
conditions at $\alp,\, \bet,\, \gam$, but many other 
choices are possible. 
\par
\section{Systems of Ordinary Differential Equations}
\par
A self-adjoint system of Sturm-Liouville operators is 
defined as an Operator $H$ acting in 
$L^{2}((\alp,\bet),\C^{m})$ according to the formula
\[
Hf_{i}(x):=-\sum_{j=1}^{m}{\rmd\over\rmd x}\left\{ a_{i,j}(x){\rmd 
f_{j}\over \rmd x}\right\} +\sum_{j=1}^{M}V_{i,j}(x)f_{j}(x).
\]
We assume that $a_{i,j}\in C^{1}[\alp,\bet]$ and $V_{i,j}\in 
L^{\infty}[\alp,\bet]$
for all $i,j$. We assume that both matrices are real 
symmetric for all $x\in [\alp,\bet]$ and that $a_{i,j}$ 
is uniformly positive definite on $[\alp,\bet]$. We finally assume 
that the operator satisfies NBC in the obvious sense for 
systems.
\par
The computation of the eigenvalues of $H$ proceeds as in the 
scalar case with one exception. Namely the removal of an 
internal NBC involves a perturbation of rank $m$ rather than 
of rank $1$ as in the scalar case. We deal with this as we did for 
higher order Sturm-Liouville operators in the last 
section. Instead of writing out the details 
in the general case, we solve a simple example, which 
exhibits the essential features of the general case.
\par
\begin{example}
Put $m=2$ and $a_{i,j}(x):=\del_{i,j}$ for all $i,j$. Let 
$\alp <0<\bet$ and let $u,v$ be arbitrary non-negative 
numbers. Then define the matrix-valued potential $V$ by
\[
V(x):=\left\{
\begin{array}{ll}
\left(\begin{array}{cc}
u&0\\0&0
\end{array}\right)&\mbox{if $\alp<x<0$}\\
\left(\begin{array}{cc}
v&v\\v&v
\end{array}\right)&\mbox{if $0<x<\bet$.}
\end{array} \right.
\]
Let $H_{1},\, H_{2},\, H$ be the operators associated with 
the above expression acting in the intervals $(\alp,0)$, 
$(0,\bet)$, $(\alp,\bet)$ respectively, all subject to NBC. 
Let $K$ be the `same' operator acting in the interval 
$(\alp,\bet)$, subject to NBC at $\alp$, $\bet$ and the 
following boundary conditions at $0$, expressed in terms of 
the operator domain: 
\begin{eqnarray*}
f_{1}(0+)&=&f_{1}(0-)\\
f_{1}^{\prime}(0+)&=&f_{1}^{\prime}(0-)\\
f_{2}^{\prime}(0+)&=& 0\\
f_{2}^{\prime}(0-)&=& 0.\\
\end{eqnarray*}
If $A_{1}$ is the operator $H_{1}+H_{2}$ acting in 
$L^{2}(\alp,\bet)$, then the quadratic forms of 
$A_{1},\, K,\, H$ 
are all given by the expression
\[
Q(f):=\int_{\alp}^{\bet}\Bigl\{ 
|f_{1}^{\prime}|^{2}+|f_{2}^{\prime}|^{2}+
\sum_{i,j=1}^{2}V_{i,j}(x)f_{i}(x)\overline{f_{j}(x)}\Bigr\}\rmd x.
\]
$A_{1}$ has the largest quadratic form domain, 
$W^{1,2}((\alp,0),\C^{2})+W^{1,2}((0,\bet),\C^{2})$, while $H$ has the 
smallest quadratic form domain $W^{1,2}((\alp,\bet),\C^{2})$, of 
codimension $2$ in the previous one. In between these 
lies the quadratic form domain of $K$, which is the set of 
$f\in  W^{1,2}((\alp,0),\C^{2})+W^{1,2}((0,\bet),\C^{2})$ such 
that $f_{1}(0+)=f_{1}(0-)$.
\par
Since each quadratic form domain is a subspace of 
codimension $1$ of the previous one, the eigenvalues of 
the operators interlace in the sense of Theorem 5. 
The eigenvalues of $H_{1},\, H_{2}$ are exactly computable, 
so these observations allow us to obtain enclosures of the 
eigenvalues of $H$ using RRTL in the standard manner.
\end{example}
\par
We have chosen this example because the eigenvalues of all 
four operators involved are essentially exactly 
computable, and it is easy to confirm the interlacing property 
directly. We put $\alp=-1$, $\bet:=2$, $u=2v=100$, and 
compute all of the eigenvalues of each operator up to 
$E:=50$.
\par
The eigenvalues of $H_{1}$ consist of all numbers of the 
form $u+n^{2}\pi^{2}/\alp^{2}$ or $m^{2}\pi^{2}/\alp^{2}$, 
where $m,\,n$ are non-negative integers. This yields the 
list:
\[
0,\, 9.870,\, 39.478,\,88.826
\]
The eigenvalues of $H_{2}$ consist of all numbers of the 
form $2v+n^{2}\pi^{2}/\bet^{2}$ or $m^{2}\pi^{2}/\bet^{2}$, 
where $m,\,n$ are non-negative integers. This yields the 
list: 
\[
0,\, 2.467,\, 9.870,\, 22.207,\, 39.478,\,61.685,\, 88.826
\]
The eigenvalues of $A_{1}$ are obtained by combining these 
two lists to obtain:
\[
0,\,0,\, 2.467,\, 9.870,\,9.870,\, 22.207,\, 39.478,\,39.478,\,
61.685,\, 88.826
\]
The eigenfunctions of $K$ and $H$ are linear combinations of 
trigonometric and exponential functions, and the eigenvalues 
are obtained by solving 
certain trancendental equations associated with the boundary 
conditions at $0$. The eigenvalues of $K$ are approximately:
\[
0,\, 0.468,\, 4.298,\, 9.870,\, 12.288,\, 24.757,\, 39.478,\, 
41.865,\, 63.639
\]
which interlace those of $A_{1}$. The 
eigenvalues of $H$ are approximately:
\[
0.449,\, 1.609,\, 4.735,\, 11.746,\,17.747,\, 27.360,\, 41.177
\] 
which interlace those of $K$.
\par
It may be seen that although the eigenvalues do interlace as 
the theory predicts, the smallest eigenvalue of $H$ is rather 
close to the second eigenvalue of $K$, a fact which does 
not help the efficiency of the TL method. The reason for 
this is that the coefficients $u,\, v$ are rather large, 
and this has the effect of partially decoupling the two 
intervals. Accurate lower bounds on the smallest eigenvalue 
of $H$ can be obtained by using a higher order version of 
the TL method.
\par

\section{Operators on graphs}
\par
The method which we have developed for Sturm-Liouville 
operators may be applied with modifications to elliptic 
partial differential operators and to discrete 
Laplacians on graphs. The first application demands the 
use of the quite complicated machinery associated with the finite 
element method. We decribe here the second 
application, which is of independent interest, and also 
involves the theory of rank $1$ 
perturbations in certain situations. 
\par
We define a graph to be a finite set $X$ together with a 
set $\cE$ of directed edges. We assume that if 
$e:=(x,y)\in\cE$ then $\overline{e}:=(y,x)\in \cE$. We define the associated 
Laplacian to be the operator acting on $l^{2}(X)$ with matrix
\[
A_{x,y}:=\left\{\begin{array}{ll}
-1&\mbox{if $(x,y)\in\cE$}\\
\degree{(x)}&\mbox{if $x=y$}\\
0&\mbox{otherwise}
\end{array}
\right.
\]
where $\degree{(x)}:=\hash\{ y\in X:(x,y)\in\cE\}$ is the degree of $x$.
\par
The quadratic form corresponding to this matrix is
\[
Q(f):={1\over 2}\sum_{(x,y)\in \cE}|f(x)-f(y)|^{2}
\]
which is a non-negative Dirichlet form, with all of the structural 
consequences of this fact. We follow standard practice in 
referring to the eigenvalues of the operator $A$ defined 
above as eigenvalues of the graph $X$.
\par
Although the matrix $A$ is finite the determination of its 
eigenvalues is not straightforward if the graph is very 
large, and one actually has the same problems in obtaining 
guaranteeed enclosures as for infinite-dimensional problems.
The first main problem is the non-existence of standard 
comparison problems. There are very few finite graphs for 
which one can compute the eigenvalues exactly, and there is 
no possibility of using a change of variables to map a 
graph to a standard soluble one as in the case of partial 
differential operators. 
\par
We present two procedures for obtaining 
enclosures of eigenvalues of finite graphs. The first applies 
to the case in which the graph is obtained from one for 
which one already has enclosures of the eigenvalues by the 
removal of a small number of chosen vertices. We leave the 
reader to formulate the corresponding lemma relating to the 
additional of a small number of vertices.
\par 
\begin{lemma} Let $Y$ be a subset of $X$ obtained by the 
removal of a small number of vertices, and let $\cG$ be 
the set of (undirected) edges of $X$ which join points of $Y$ to points 
of $X\backslash Y$. Then one may compute enclosures of the eigenvalues of 
$Y$ from those of $X$ in $\hash (\cG)$ homotopy steps.  
\end{lemma}
\par
\Proof Let $\{ e_{i}\}_{i=1}^{n}$ be some enumeration of the 
edges in $\cG$. Let $A_{i}$ be the matrix acting in 
$l^{2}(X)$ which is obtained from that of $A$ by the removal of the 
contribution of the edges $e_{1},\ldots, e_{i}$. Then
\[
A\geq A_{1}\geq \ldots \geq A_{n}
\]
and each matrix is a rank $1$ perturbation of the next one 
in the chain. It follows by an argument similar to that of 
Theorems 1 and 5 that the eigenvalues of $A_{i}$ interlace 
those of $A_{i+1}$ for every $i$. This is the fundamental 
requirement for transferring accurate enclosures of the eigenvalues 
from $A_{i}$ to $A_{i+1}$ by the RRTL method. The operator 
$A_{n}$ is the direct sum of the discrete Laplacians of $Y$ 
and $X\backslash Y$. Since we have assumed that $X\backslash 
Y$ is small, its eigenvalues may be computed independently 
by a direct procedure. Removing these eigenvalues leaves 
those of $Y$.
\par  
Let $X$ be a 
finite subset of $\Z^{N}$ for some $N$. The edges of $X$ 
are defined to be those pairs $x,y\in X$ such that 
\[
\sum_{r=1}^{N}|x_{r}-y_{r}|=1.
\] 
In this situation we have the bound $\degree (x)\leq 2N$ for 
all $x\in X$. The operator $A$ may be considered to be the 
discrete Laplacian on $X$ subject to Neumann boundary 
conditions, since $0$ is always an eigenvalue of $A$, the 
corresponding eigenfunction being constant. The multiplicity 
of the eigenvalue $0$ equals the number of connected 
components of the graph.
\par
\begin{example}
Let $X\subset \Z^{2}$ be the set
\[
\{ (m,n):1\leq m\leq k, 1\leq n\leq k\}
\]
and let $Y$ be obtained by the removal of the set
\[
Z:=\{ (1,1),(1,2),(2,1)\}.
\] 
Then an enumeration of the (undirected) edges of 
$\cG$ is $e_{1}:=((1,2),(2,2))$, 
$e_{2}:=((2,1),(2,2))$, $e_{3}:=((1,2),(1,3))$, 
$e_{4}:=((2,1),(3,1))$. We chose $k:=7$ and computed the 
smallest $6$ eigenvalues of the operators $A_{i}$.
The interlacing property is verified. The eigenvalues of 
$A_{4}$ coincide with those of the discrete Laplacian of 
$Y$, together with the eigenvalues 
$0,1,3$ of the Laplacian of $Z$. The numbers in the table 
below are actually $k^{2}$ times the 
eigenvalues, so that they may be compared with the 
eigenvalues of $-\lap$ on the unit square subject to NBC. 
\[
\begin{array}{ccccccc}
&\mu_{0}&\mu_{1}&\mu_{2}&\mu_{3}&\mu_{4}&\mu_{5}\\
A&0&9.705&9.705&19.410&36.898&36.898\\
A_{1}&0&9.515&9.705&19.142&33.499&36.898\\
A_{2}&0&9.361&9.574&18.367&30.187&34.782\\
A_{3}&0&5.868&9.540&13.119&23.836&34.571\\
A_{4}&0&0&9.095&11.471&23.049&32.525
\end{array}
\]
The interlacing property states that the number $\lam$ immediately 
above any eigenvalue $\mu$ of $A_{i}$ in the table is a 
lower bound for the next eigenvalue of $A_{i}$.
Let us suppose that the values in the first row of the table 
are close to accurate enclosures of the eigenvalues of $A$, and 
that all of the other entries in the table
are rigorous upper bounds, which we  
expect to be accurate. Using the interlacing property we 
deduce that $\mu_{5}(A_{1})=36.898$. The fact that the
(upper bounds on the) 
other eigenvalues of $A_{1}$ are widely separated enables 
us to use TL to confirm that they 
have been found accurately. Interlacing establishes that 
$\mu_{5}(A_{2})\geq 33.499$, and we confirm that the other 
eigenvalues of $A_{2}$ have been computed accurately as 
before. If interval arithmetic has been used we end up with 
accurate enclosures of all of the eigenvalues of $A_{4}$ up 
to and including $\mu_{4}$. The final reason why this 
procedure works is that the entries in the column 
labelled $\mu_{4}$ and the row labelled $A_{4}$ 
decrease rapidly enough for TL to be 
an efficient method.  
\par  
The above example is purely illustrative: the same procedure 
can be carried out for 
values of $k$ which are large enough for the problem of 
obtaining eigenvalue enclosures to be non-trivial. If at 
some stage an eigenvalue does not decrease enough from one 
stage to the next for TL, 
then either we have to proceed to a lower eigenvalue, or we 
must use a higher order version of TL.
\end{example}
\par
The above method is not suitable for obtaining enclosures 
of the eigenvalues of a graph which is far from any graph for which 
eigenvalue enclosures are already known. As a typical 
example we mention the set of all $(m,n)\in \Z^{2}$ which 
satisfy all three inequalities $m^{2}+n^{2}<16d^{2}$, 
$(m-2d)^{2}+n^{2}>d^{2}$ and $(m+2d)^{2}+n^{2}>d^{2}$ 
where $d$ is a large positive number. 
\par
In cases such as the above we 
combine the continuous homotopy procedure  
introduced by Goerisch and Plum with the hierarchical 
homotopy method we introduced for Sturm-Liouville operators. 
The idea is to subdivide 
$X$ into several more or less convex parts each of which is small enough that 
eigenvalue enclosures can be obtained by a direct method. 
These parts are then joined together in pairs as described below, 
obtaining eigenvalue enclosures for 
the larger parts. If the initial subdivision is into $k:=2^{N}$ 
parts, then after the first stage one has $2^{N-1}$ parts, 
and the procedure terminates after $N$ stages. It remains 
to describe how to join together two subgraphs.
\par
Let $X=Y\cup Z$ where $Y,\, Z$ are disjoint subgraphs,  
let $\cG$ be the set of edges joining points of $Y$ and 
$Z$, and let $\cF$ be the complement of $\cG$ in the set 
$\cE$ of all edges of $X$. Given $0\leq s\leq 1$, let 
$A_{s}$ be the matrix associated with the quadratic form
\[
Q_{s}(f):={1\over 2}\sum_{(x,y)\in \cF}|f(x)-f(y)|^{2} 
+{s\over 2} \sum_{(x,y)\in 
\cG}|f(x)-f(y)|^{2}. 
\] 
\par
\begin{lemma}
The eigenvalues of $A_{s}$ are increasing real analytic 
functions of the parameter $s$.  The eigenvalue list of 
$A_{0}$ is just the union of the two eigenvalue lists of $Y$ and $Z$. 
At the other end $A_{1}$ is the discrete Laplacian of $X$.
\end{lemma} 
\par
\Proof The first statement is part of received knowledge 
\cite{Ka}, while the second depends upon the fact 
that  $A_{0}$ is the direct sum of the discrete 
Laplacians of $Y$ and $Z$.
\par
The procedure for obtaining eigenvalue enclosures for $X$ is 
similar to that of Goerisch and Plum \cite{G,P1}. We consider the operators 
$A_{s(r)}$ for a large enough chosen sequence 
$0=s_{0}<s_{1}<\ldots <s_{p}=1$. If we have enclosures of the 
eigenvalues of $A_{s(i)}$ then these provide lower bounds on 
the eigenvalues of $A_{s(i+1)}$ which may be adequate to 
obtain enclosures of the eigenvalues of the latter 
operator by the RRTL procedure. Eigenvalue crossings may 
occur at certain values of $s$, but these are handled using 
the higher order TL procedure.
\par
\begin{example}
Let $X:=Y\cup Z\subset \Z^{2}$ where
\begin{eqnarray*}
Y&:=&\{ (m,n):1\leq x \leq h-1,\, 1\leq y\leq x \}\\
Z&:=&\{ (m,n):h\leq x\leq 2h-1,\, 1\leq y\leq 2h-x\}
\end{eqnarray*}
where $h$ is some positive integer. Let the set $\cE$ of 
edges of $X$ be those inherited from $\Z^{2}$ as before. 
The undirected edges of $\cG$ are 
of the form $(h-1,r),(h,r)$ where $1\leq r \leq h-1$, and 
separate the triangle $X$ into two smaller triangles. We 
list the $7$ smallest eigenvalues of $A_{s}$ below for $h:=8$ and 
$s:=0,\, 0.2,\, 1$. A larger number of values of $s$ were originally 
computed, but these are the only ones needed.
\[
\begin{array}{cccccccc}
&\mu_{0}&\mu_{1}&\mu_{2}&\mu_{3}&\mu_{4}&\mu_{5}&\mu_{6}\\
A_{0}&0&0&0.12061&0.15224&0.25330&0.32139&0.46791\\
A_{0.2}&0&0.04705&0.13054&0.23249&0.27273&0.42485&0.49950\\
A_{1}&0&0.07244&0.13259&0.27719&0.33076&0.51058&0.60389
\end{array}
\]
Each eigenvalue $\mu_{i}$ is a monotonic increasing function 
of $s$. Suppose that we already know that the numbers in 
the first row are accurate approximations to the eigenvalues 
of $A_{0}$, and that the other numbers are rigorous upper 
bounds to the corresponding eigenvalues, as determined by 
RR. We use the monotonicity to deduce that 
$\mu_{6}(A_{0.2})\geq 0.46791$. Using TL we then confirm 
that $\mu_{5}(A_{0.2})$ is accurate. Monotonicity implies 
that $\mu_{5}(A_{1})\geq 0.42485$ and we are finally able to 
confirm the accuracy of $\mu_{j}(A_{1})$ for $j=4,3,2,1,0$ 
in turn by TL. If all of the computations have been done in 
interval arithmetic, we have obtained enclosures of the 
eigenvalues of $A_{1}$ from those of $A_{0}$.
\par
The values of $\mu_{1}(A_{s})$ for $s=0..1(0.1)$ are 
\begin{eqnarray*}
&0,\,& 0.03153,\, 0.04705,\, 0.05559,\, 0.06085,\, 0.06439\\ 
&&0.06693,\, 0.06882,\, 0.07030,\, 0.07148,\, 0.07244
\end{eqnarray*}
This list exhibits a common pattern of rapid increase for 
small values of $s$ followed by little change for large values. 
In fact it follows from RR that $\mu_{1}(A_{s})$ is a 
concave function of $s$, but this need not be true for higher 
eigenvalues. 
\end{example}
\par
A precondition for applying the above method is that eigenvalue enclosures 
of the parts $X_{1},\ldots ,X_{k}$ into which we subdivide $X$ 
should already be known. This may be achieved by making each 
$X_{i}$ small enough so that all of its eigenvalues can  be 
computed by a direct method. If some of the parts are 
rectangles or one of a very small number of other graphs then 
their eigenvalues may be exactly known. The following 
argument shows that for some purposes we may dispense with 
knowledge of accurate enclosures of the eigenvalues of the parts entirely. 
It requires instead an assumption about the geometry of the 
parts, expressed initially in terms of a lower bound on their first 
non-zero eigenvalues.
\par
Given $a>0$ we say that a finite graph $(X,\cE)$ lies in 
$\cC_{a}$ if its first non-zero eigenvalue $\mu_{1}$ 
satisfies
\[
\mu_{1}\geq a\, d(X,\cE)^{-2}
\]
where $d(X,\cE)$ is the diameter of the graph. We investigate the geometric 
significance of this condition in the next section. The 
value of $b$ in the following theorem indicates how small the 
individual subsets in a partition of $X$ need to be in order to 
be able to obtain accurate enclosures of the eigenvalues of 
$X$ without already possessing accurate enclosures of the 
eigenvalues of the subsets.
\par
\begin{theorem}
Suppose that $a>0$ and that $\{ X_{i} \}_{i=1}^{k}$ is a partition of the 
graph $X$, each subset of which lies in $\cC_{a}$. Suppose also 
that $b>0$ and that  
\[
d(X_{i},\cE_{i})\leq {1\over b}d(X,\cE)
\]
for all $1\leq i \leq k$. Then one may obtain accurate 
enclosures of all eigenvalues of $X$ less than 
\[
E:=ab^{2}\, d(X,\cE)^{-2}
\]
by a continuous homotopy method.
\end{theorem}
\par
\Proof The conditions of the theorem imply that the first 
non-zero eigenvalue of each part is at least as big as 
$E$. Let $(Y,\cF)$ be the union of two of 
the subgraphs $(X_{i},\cE_{i})$, and let 
$(Y,\overline{\cF})$ be obtained by including also those 
edges of the graph $(X,\cE)$ which connect the two parts of 
$Y$. The 
spectrum of $(X,\cF)$ below $E$ consists of the 
eigenvalue $0$ with multiplicity $2$ and nothing else. This 
precise if unusual information is enough to obtain accurate 
enclosures of the 
spectrum of $(Y,\overline{\cF})$ below $E$ by the Goerisch-Plum continuous homotopy 
method. By a repetition of this method, carried out in a 
hierarchical manner, one eventually obtains accurate 
enclosures of the spectrum of $(X,\cE)$ below $E$.
\par
\section{The Poincar\'e Inequality}
\par
In order to implement the above ideas one needs to obtain 
geometric conditions which imply 
that the first non-zero eigenvalue of a connected graph $(X,\cE)$ 
has a lower bound 
of order $d^{-2}$, where $d$ is the diameter of the graph. 
Examples in \cite{DS} show that this is not always uniformly true for a 
family of graphs parametrised by the diameter as 
$d\to\infty$, 
but their discrete version of the Poincar\'e inequality can 
be rewritten to provide exactly what we need. We develop the 
theory of this section at a greater level of generality 
than before, because of its independent interest.
\par
Let $b:\cE\to (0,\infty)$ be a positive weight function 
satisfying $b(e)=b(\overline{e})$ for all $e\in\cE$. Let 
$|X|$ denote the number of points in $X$. Given a path 
$\gam:=(\gam_{1},\ldots,\gam_{k})$, we define its length 
to be 
\[
|\gam |:=\sum_{i=1}^{k}b(\gam_{i-1},\gam_{i})^{-1}
\] 
and then define the diameter $d$ of $X$ using this notion of 
length, in the usual manner.
\par
Let $A$ 
be the operator on $l^{2}(X)$ associated with the quadratic 
form
\[
Q(f):={1\over 2}\sum_{(x,y)\in \cE}b(x,y)|f(x)-f(y)|^{2}.
\]
It is easily seen that $0$ is an eigenvalue of $A$ of 
multiplicity $1$, the corresponding normalised eigenfunction 
satisfying 
$\phi_{0}(x)=|X|^{-1/2}$ for all $x\in X$.
\par
In order to obtain a lower bound on the first non-zero 
eigenvalue $\mu_{1}$ of $A$, we follow closely Diaconis and Stroock 
\cite{DS} (and Poincar\'e). Suppose that $\Gam$ is a set of paths in $X$, one 
path from $\Gam$ joining every ordered pair of points $x,y\in X$. We impose two 
constraints on the choice of this set of paths. The first, that 
\[
|\gam_{x,y}|\leq \alp d
\]
for some $\alp$ and all $x,y\in X$, is self-explanatory. 
The second is that 
\[
\hash\{ \gam\in \Gam: e\in\gam\}\leq \bet d |X|
\]
for some $\bet >0$ and all $e\in \cE$. Since the total number of 
paths in $\Gam$ is $|X|^{2}$, this is a constraint on how well 
distributed the paths are. The assumption is in precisely 
the form needed for applications.
\par
\begin{theorem}
Under the above two assumptions we have
\[
\mu_{1} \geq {1\over \alp\bet d^{2}}.
\]
\end{theorem}
\Proof  If $e:=(x,y)\in\cE$ we put $\partial f(e):=f(y)-f(x)$. We have
\begin{eqnarray*}
|X|\,\norm \phi -\la \phi,\phi_{0}\ra 
\phi_{0}\norm^{2}
&=&{1\over 2}\sum_{x,y\in X}|\phi(x)-\phi(y)|^{2}\\
&=&{1\over 2}\sum_{x,y\in X}|\sum_{e\in\gam_{x,y}}\partial\phi(e)|^{2}\\
&\leq&{1\over 2}\sum_{x,y\in 
X}|\gam_{x,y}|\sum_{e\in\gam_{x,y}}b(e)|\partial\phi(e)|^{2}\\
&\leq&KQ(\phi)
\end{eqnarray*}
where
\begin{eqnarray*}
K:&=&\sup_{e\in\cE}\sum_{\gam_{x,y}\ni e}|\gam_{x,y}|\\
&\leq&\alp d\hash\{  \gam\in \Gam: e\in\gam   \}\\
&\leq&\alp\bet d^{2}|X|.
\end{eqnarray*}
The proof is completed by using the variational 
characterisation of $\mu_{1}$:
\[
\mu_{1}=\inf\left\{ 
{Q({\phi})  \over   \norm \phi -\la \phi,\phi_{0}\ra \phi_{0}\norm^{2}}
: 0\not=\phi\in l^{2}(X) \right\}.
\]
\par
Upper bounds of a similar type on $\mu_{1}$ are relatively easy to obtain by 
applying the variational inequality to suitable 
test functions, but we do not need them here.
\par
The following application of the above theorem is a typical 
building block for the implementation of the ideas in the 
last section.
\begin{theorem}
Define $X\subset \Z^{2}$ by
\[
X:=\{ (i,j):1\leq i \leq n ,\,\, 1\leq j\leq f(i)\}
\]
where $f\{ 1,\ldots,n\} \to \{ 1,\ldots m\}$ is a non-decreasing 
function. Let $A$ be the discrete Laplacian on $X$, 
corresponding to the choice $b\equiv 1$ above. Then 
\[
\mu_{1}\geq {1\over (n+f(n)-2)\max\{ n,f(n) \}} \geq d^{-2}
\]
where $d$ is the diameter of $X$. 
\end{theorem}
\par
\Proof It follows from the definition of $X$ that $d=n+f(n)-2$ 
and $n\leq |X| \leq nf(n)$. Our main task 
is to define the set $\Gam$ of paths. If $i\leq i^{\prime}$ 
then the path from $(i,j)$ to $(i^{\prime},j^{\prime})$ is 
the horizontal line from  $(i,j)$ to $(i^{\prime},j)$, 
followed by the vertical line from  $(i^{\prime},j)$ to 
$(i^{\prime},j^{\prime})$. Because $f$ is monotonic, this 
is entirely contained in $X$. If $i\geq i^{\prime}$ we take 
a similar path. It may be seen that every path $\gam$ has 
length at most $n+f(n)-2$. The number of paths 
through any horizontal edge $e\in\cE$ is at most $n|X|$, 
while the number through any vertical edge is at most 
$f(n)|X|$. A slight modification of the estimate of $K$ in 
the last theorem completes the proof.
\par
The number of paths through any edge 
$e\in \cE$ can be bounded more efficiently if further 
information about $f$ is provided. If $f(i):=1$ for $1\leq 
i\leq n-1$ and $f(n):=n$, then 
the theorem
yields $\mu_{1}\geq {1\over n(2n-2)}$ while one 
actually has $\mu_{1}\sim {\pi^{2}\over 4n^{2}}$ for large 
$n$. It would be valuable to determine the largest 
constant $c$ such that $\mu_{1}\geq c/n^{2}$ for all graphs 
of the type described in the above theorem, subject to 
$f(n)\leq n$. It appears that even the continuous analogue 
of this problem is unsolved.   
\section{The Laplacian in $N$ Dimensions}
\par
We finally describe the modifications to the above ideas 
needed to provide enclosures of the eigenvalues of a partial 
differential operator. We will consider only the case in 
which $H:=\lap$, acting in $L^{2}(\Ome)$ subject to NBC, 
where $\Ome$ is a bounded region in $R^{N}$ with piecewise 
smooth boundary. However, the same method applies to 
variable coefficient elliptic operators subject to other boundary 
conditions. By the eigenvalues of any region $\Ome$ 
we mean the eigenvalues of $-\lap$ acting in $L^{2}(\Ome)$ 
subject to NBC.
\par
If the region $\Ome$ is diffeomorphic to the unit 
ball $B$ 
Plum \cite{P2} has described a method of obtaining enclosures 
on the eigenvalues by transferring the operator to 
$L^{2}(B,m(x)\rmd x)$ where $m$ is a suitable positive 
weight, and then using his coefficient homotopy method. 
This cannot be adapted to treat the case in which $\Ome$ 
contains more than one hole. The method described below is capable of dealing 
with regions containing any number of holes.
\par
Suppose we wish to find all of the eigenvalues of a region $\Ome$ 
which are smaller than a given number $E>0$. The first step is to 
divide the region into subregions $\{ \Ome_{i} \}_{i=1}^{k}$ each with 
a piecewise smooth boundary.  We assume that enclosures of the  
eigenvalues less than $E$ of each subregion are known, either because its 
eigenvalues are exactly computable or because it 
is small enough with a regular enough shape for $0$ to be 
its only eigenvalue below $E$. We also assume that the 
subregions can be recombined in pairs in a hierarchical 
fashion to recover the original set $\Ome$. The task 
therefore is to obtain enclosures of the eigenvalues of the union of two 
regions which have some common boundary when we 
already have enclosures of the eigenvalues of the individual regions.
\par
Let $U$, $V$ be disjoint bounded connected regions in $\R^{N}$ with 
piecewise smooth boundaries and suppose that their common 
boundary $B$ is a non-empty $(N-1)$-dimensional surface. 
Put $\Ome:=U\cup V\cup B$. Suppose also that we have 
accurate enclosures of all the eigenvalues of each region up 
to the number $E$. If we combine the two lists of 
eigenvalues into a single increasing list $\{ \mu_{i}(0)\}$, 
then this list provides the eigenvalues of the operator 
$H_{0}=-\lap$ acting in $L^{2}(\Ome)$ subject to NBC on 
$\partial U\cup\partial V$. Note that the number $0$ is an 
eigenvalue of multiplicity $2$.
\par
We introduce a family of quadratic forms $Q_{s}$ defined 
for $0\leq s <\infty$, all having the same quadratic form 
domain $\cD_{0}:=W^{1,2}(U)+W^{1,2}(V)$. A core for this 
subspace consists of all functions on $\overline \Ome$ which are $C^{1}$ 
except that they are allowed to be discontinuous as one crosses 
$B$. Every function $f\in \cD_{0}$ has $L^{2}$ boundary values 
on $B$ which we denote by $f_{\pm}$ depending upon which 
side of $B$ one approaches it from. We then define
\[
Q_{s}(f):=\int_{\Ome}|\nabla f|^{2} +s\int_{B}|f_{+}-f_{-}|^{2}
\]
where the second integral is with respect to the natural 
surface measure on $B$. It is evident that the forms $Q_{s}$ 
are monotonic increasing. It may be shown that the 
perturbation term is relatively compact, so that the right-hand side is 
the closed form associated with a certain non-negative 
self-adjoint operator $H_{s}$. The eigenvalues of $H_{s}$ 
are increasing real-analytic functions of $s$ by \cite{Ka}.
\par
Since the quadratic forms are monotonic increasing as a 
function of $s$ they converge to a limit $Q$ defined by
\[
Q(f):=\lim_{s\to+\infty}Q_{s}(f)
\]
where we adopt the standard convention \cite{D1} that $Q(f)=+\infty$ 
whenever $f$ does not lie in the form domain $\cD$ of $Q$. It is 
clear that 
\[
\cD=\{ f\in \cD_{0}: f_{+}=f_{-}\mbox{ on } B\}.
\]
This space is exactly $W^{1,2}(\Ome)$, so the operator 
associated with $Q$ is $H:=-\lap$ acting in 
$L^{2}(\Ome)$ subject to NBC on $\partial\Ome$. An immediate 
consequence is that \[
\lim_{s\to +\infty}\mu_{n}(s)=\mu_{n}
\]
for every $n$, the limit being monotone.
\par
Having chosen a suitable increasing sequence 
$0=s_{0}<s_{1}<\ldots<s_{n}$ one obtains accurate 
enclosures of the eigenvalues of each $H_{s_{i+1}}$ from those of $H_{s_{i}}$ 
by the Goerisch-Plum homotopy method. If $s_{n}$ is large 
enough then the eigenvalues of $H_{s_{n}}$ will be good enough 
lower bounds of the eigenvalues of $H$ to enable us to apply 
RRTL to obtain accurate enclosures of the eigenvalues of $H$. 
Many of the details are similar 
to what we have already discussed, and we concentrate on 
the novelties. 
\par 
The upper bounds on all eigenvalues 
are obtained by RR using a suitable test function space lying 
in the quadratic form domain $\cD_{0}$ or $\cD$. An 
obvious choice is to use a finite 
element subspace in which the elements are linear, but more 
sophisticated elements are probably needed for accurate 
results. The 
continuity requirement for two elements which have a common edge 
is suspended if that edge lies within $B$. The restriction 
of the quadratic form $Q_{s}$ to the test function space 
can be expanded in terms of the values of the 
elements at the vertices, noting that there will be two 
values at each vertex lying on $B$, one corresponding to 
each side of $B$. 
\par
The required lower bounds can only be obtained by TL if one 
takes a test function space lying in the operator 
domain, which is different for each operator, even though 
every operator $H_{s}$ is equal to $-\lap$ on its own 
domain. One can use a finite element subspace consisting of $C^{2}$ 
functions, but this has to respect 
not only the Neumann boundary condition on $\partial W$ but 
also certain $s$-dependent internal boundary conditions on $B$.
\par
Let $\partial f_{\pm}$ denote the normal derivatives of 
$f$ on the two sides of $B$, both taken in the direction 
from the $-$ side of $B$ to the $+$ side.
An application of Gauss' theorem shows that 
the internal boundary condition is 
\[
\partial f_{+ }(x) =\partial f_{-}(x)=s\{ f_{+}(x)-f_{-}(x)\}
\] 
for all $x\in B$.
\par
It is not easy to demonstrate that the above theory works well in practice, 
without a substantial amount of effort writing the relevant 
code. An appropriate choice of the test function space is 
crucial if one is to get good enclosures, as may be the use 
of a preconditioning procedure, and we leave this 
to a future publication.    
\vskip 0.3in
\par
{\bf Acknowledgments } We would like to thank M Plum for 
several valuable exchanges in the course of this work.
\vskip 0.3in

E.Brian.Davies@kcl.ac.uk\newline
King's College, London

\begin{thebibliography}{99}
\bibitem{D1} E B Davies. Spectral Theory and Differential 
Operators. Cambridge Univ. Press, 1995.
\bibitem{D2} E B Davies. Spectral enclosures and complex 
resonances for general self-adjoint operators. Preprint 1997.
\bibitem{DS} P Diaconis and D W Stroock: Geometric bounds 
for eigenvalues of Markov chains. Ann. Appl. Prob. 1 (1991) 
36-61.
\bibitem{G} F Goerisch: Ein Stufenverfahren zur Berechnung 
von Eigenwertschranken. In `Numerical Treatment of 
Eigenvalue Problems' vol. 4 ISNM 83 (Ed. J Albrecht et al) 
pp 104-114. Birkhauser-Verlag, Basel, 1987.
\bibitem{Ka} T Kato: Perturbation Theory of Linear Operators. 
Springer-Verlag, Berlin, Heidelberg, New York, 1966. 
\bibitem{Lo} R Lohner: Verified solution of eigenvalue 
problems in ordinary differential equations. Unpublished 
manuscript, 1990.
\bibitem{P1} M Plum: Eigenvalue inclusions for second order 
ordinary differential operators by a numerical homotopy 
method. J. Appl. Math. and Phys. 41 (1990) 205-226.
\bibitem{P2} M Plum: Bounds for eigenvalues of second order 
elliptic operators. J. Appl. Math. and Phys. 42 (1991) 
848-863. 
\bibitem{P3} M Plum: Guaranteed numerical bounds for 
eigenvalues. In `Spectral Theory and Computational Methods 
of Sturm-Liouville Problems', eds. D Hinton and P W Schaefer. 
Marcel Dekker, New York, Basel, 1997.
\bibitem{ZM}  S Zimmerman and U Mertins: Variational bounds to 
eigenvalues of self-adjoint eigenvalue problems with 
arbitrary spectrum. Zeit. f\"ur Anal. und ihre Anwendungen. 
J. for Analysis and its Applications 14 (1995) 327-345.

\end{thebibliography}
\end{document}